\documentclass[11pt,fceqn]{article}
\usepackage{amsmath}
\usepackage{amsfonts}
\usepackage{cite}
\usepackage{amssymb,graphics,graphicx,subfigure,caption2,rotating}
\usepackage{amsmath, latexsym}
\usepackage[amsmath,thmmarks]{ntheorem}
\numberwithin{equation}{section}
\newtheorem{theorem}{Theorem}[section]
\newtheorem{lemma}{Lemma}[section]

\newcommand{\beq}{\begin{equation}}
\newcommand{\eeq}{\end{equation}}
\newcommand{\beqn}{\begin{eqnarray}}
\newcommand{\eeqn}{\end{eqnarray}}
\date{}
\begin{document}

\date{}
\title{A class of maximum-based iteration methods for the generalized absolute value equation\thanks{This research was subsidized
by National Natural Science Foundation of China
(No.11961082,12131004).}}
\author{Shiliang Wu\thanks{Corresponding author:
slwuynnu@126.com},  Deren Han\thanks{handr@buaa.edu.cn}, Cuixia
Li\thanks{lixiatk@126.com}\\
{\small{\it  $^{\dagger,\S}$School of Mathematics, Yunnan Normal University, }}\\
{\small{\it  Kunming, Yunnan, 650500, P.R. China}}\\
{\small{\it $^{\dagger}$Yunnan Key Laboratory of Modern Analytical Mathematics and Applications,}}\\
{\small{\it  Yunnan Normal University, Kunming, Yunnan, 650500, P.R.
China}}
\\
{\small{\it $^{\ddagger}$LMIB of the Ministry of Education, School of Mathematical Sciences,}}\\
{\small{\it  Beihang University, Beijing, 100191, P.R. China}} }
\maketitle
\begin{abstract}
In this paper, by using $|x|=2\max\{0,x\}-x$, a class of
maximum-based iteration methods is established to solve the
generalized absolute value equation $Ax-B|x|=b$. Some convergence
conditions of the proposed method are presented. By some numerical
experiments, the effectiveness and feasibility of the proposed
method are confirmed.


\end{abstract}

\section{Introduction}
The generalized absolute value equation (GAVE) is to look for $x\in
\mathbb{R}^{n}$ such that
\begin{equation}\label{eq:11}
Ax-B|x|=b \ \mbox{with}\ A,B\in \mathbb{R}^{n\times n} \mbox{and}\
b\in \mathbb{R}^{n},
\end{equation}
where $|\cdot|$ denotes the absolute value. The GAVE (\ref{eq:11})
often arises in diverse fields, like  flows in porous media, linear
programming, nonnegative constrained least squares problems,
quadratic programming, complementarity problem, bimatrix games, and
others, see \cite{Rohn,Mangasarian,Mangasarian2,Wu,Cottle,Zheng}.

Nowdays, for the theoretical analysis of the GAVE (\ref{eq:11}),
there exist many interesting results. For example, Wu and Shen in
\cite{Wu5} presented some necessary and sufficient conditions for
the unique solvability of the GAVE (\ref{eq:11}). Wu and Li
\cite{Wu6} discussed the error and perturbation bounds of the GAVE
(\ref{eq:11}).  Hlad\'{\i}k in \cite{Hla23} investigated the
properties of the solution set of the GAVE (\ref{eq:11}) with $B=I$
(the didentity matrix). Zamani and Hlad\'{\i}k in \cite{Zamani23}
discussed its error bounds and a condition number.

Besides that, to efficiently solve the GAVE (\ref{eq:11}), many
efficient iteration methods have been designed, including the
successive linearization method \cite{Mangasarian}, the sign accord
method \cite{Rohn2}, the Picard and Picard-HSS methods
\cite{Salkuyeh},  the generalized Newton (GN) method
\cite{Mangasarian3}, and so on. Recently, Wang et al.  in
\cite{Wang} designed the modified Newton (MN) method, see  Method 1.

\textbf{Method 1} Assume that $x^{0}\in \mathbb{R}^{n}$ is an
arbitrary initial guess. For $k = 0,1,2,\ldots$ until the iteration
sequence $\{x^{k}\}_{k=0}^{\infty}\in\mathbb{R}^{n}$ is convergent,
computing $x^{(k+1)}$ by
\begin{equation}\label{eq:12}
x^{k+1}=(A+\Omega)^{-1}(\Omega x^{k}+B|x^{k}|+b), k=0,1,,2,\ldots,
\end{equation}
where $A+\Omega$ is invertible and $\Omega$ is positive
semi-definite.

To expedite the convergence speed of Method 1, by making use of the
matrix splitting technique, i.e., by expressing $A$ as $A=M-N$, Zhou
\emph{et al.} in \cite{Zhou21} presented the following Newton-based
matrix splitting (NMS) method, see Method 2.

\textbf{Method 2}  Assume that $x^{0}\in \mathbb{R}^{n}$ is an
arbitrary initial guess. For $k = 0,1,2,\ldots$ until the iteration
sequence $\{x^{k}\}_{k=0}^{\infty}\in\mathbb{R}^{n}$ is convergent,
computing $x^{k+1}$ by
\begin{equation} \label{eq:13}
\begin{split}
x^{k+1}=(M+\Omega)^{-1}((N+\Omega)x^{k}+B|x^{k}|+b), k=0,1,2\ldots,
 \end{split}
 \end{equation}
where $\Omega+M$ is invertible and $\Omega$ is a given matrix.

Comparing Method 2 with Method 1, the former not only covers the
latter, but also may be performed easily by choosing the suitable
matrix $\Omega$ and matrix splitting of $A$, moreover, may be fit
for the indefinite matrix $A$ as well.

In this paper, we continuously develop the numerical method for
solving the GAVE (\ref{eq:11}) on condition that the GAVE
(\ref{eq:11}) has a unique solution (see \cite{Wu5}). By using
$|x|=2\max\{0,x\}-x$ for the GAVE (\ref{eq:11}), a class of
maximum-based iteration methods is established, and is completely
different from the above existing numerical methods. Not only that,
the building new iteration method does not contain $|x|$, and takes
full advantage of information of matrix $B$. In particular, under
certain condition, determining the solution of the absolute value
equation is equal to determining the solution of the linear systems.
We investigate the convergence properties of the proposed method.
Further, by using some numerical experiments, the effectiveness and
feasibility of the proposed method are confirmed.

The remainder of the paper is organized below.  In Section 2, a
class of maximum-based iteration methods is established, and its
convergence properties are discussed. In Section 3, the numerical
results for the proposed method are reported. In Section 4, we draw
some conclusions to finish this paper.

\section{A class of maximum-based iteration methods}

Firstly, for the later discussion, we recall some necessary
terminology and lemmas.

For $A=(a_{ij}) \in \mathbb{R}^{n\times n}$, it is called as a
$Z$-matrix if $a_{ij}\leq0$ for $i\neq j$;  a nonsingular $M$-matrix
if $A^{-1}\geq0$ and  $A$ is a $Z$-matrix;  an $H$-matrix if its
comparison matrix $\langle A\rangle=(\langle a_{ij}\rangle) \in
\mathbb{R}^{n\times n}$ ($\langle a_{ii}\rangle=|a_{ii}|$ and
$\langle a_{ij}\rangle=-|a_{ij}|$ for $i\neq j$) is  a nonsingular
$M$-matrix; an $H_{+}$-matrix if it is an $H$-matrix with positive
diagonal entries, see\cite{Bai}. $\rho(A)$ and $\|A\|_{2}$ in order
denote the spectral radius and the 2-norm of the matrix $A$.

Let $A = M- N$ with $\mbox{det}(M)\neq 0$. Then it is named an
$M$-splitting if $N \geq 0$  and $M$ is a nonsingular $M$-matrix.
Further, if  $A = M - N$ is an $M$-splitting with $A$ being a
nonsingular $M$-matrix, then $\rho(M^{-1}N) < 1$, see \cite{Berman}.

\begin{lemma} \emph{\cite{Frommer}}
Let $A = D-B$, where $A\in \mathbb{R}^{n\times n}$ is an $H$-matrix
and  $D$ is its diagonal part. Then matrices $A$ and $|D|$ are
nonsingular,  $|A^{-1}|\leq\langle A\rangle^{-1}$.
\end{lemma}

Secondly, we will establish a class of maximum-based iteration
methods to solve the GAVE (\ref{eq:11}). To this end, by using
$|x|=2\max\{0,x\}-x$ for the GAVE (\ref{eq:11}), the GAVE
(\ref{eq:11}) is equal to
\begin{equation}\label{eq:21}
(A+B+\Omega)x=\Omega x+2B\max\{0,x\}+b,
\end{equation}
where $A+B+\Omega$ is invertible. Hence, based on Eq. (\ref{eq:21}),
a class of maximum-based iteration methods can be established, see
Method 3.

\textbf{Method 3}  Assume that $x^{0}\in \mathbb{R}^{n}$ is an
arbitrary initial guess. For $k = 0,1,2,\ldots$ until the iteration
sequence $\{x^{k}\}_{k=0}^{\infty}\in\mathbb{R}^{n}$ is convergent,
computing $x^{k+1}$ by
\begin{equation} \label{eq:22}
\begin{split}
x^{k+1}=(A+B+\Omega)^{-1}(\Omega x^{k}+2B\max\{0,x^{k}\}+b),
 \end{split}
 \end{equation}
where matrix $A+B+\Omega$ is invertible.

Clearly, Method 3 gives a new general framework for solving the GAVE
(\ref{eq:11}). Similar to Method 2, some relaxation methods are
generated by selecting the matrix splitting of matrix $A+B$ for
(\ref{eq:21}). Here is omitted. It is easy to see that Method 3 does
not contain $|x|$ and takes full advantage of information of matrix
$B$.

In addition, it is noted that if $x$ is non positive, i.e., $x\leq
0$, then Eq. (\ref{eq:21}) reduces to the linear systems
\begin{equation*}
(A+B+\Omega)x=\Omega x+b\ \mbox{or} \ (A+B)x=b;
\end{equation*}
if $x$ is positive, i.e., $x>0$, then Eq. (\ref{eq:21}) reduces to
the linear systems
\begin{equation*}
(A-B+\Omega)x=\Omega x+b\ \mbox{or} \ (A-B)x=b.
\end{equation*}
This shows that under certain condition, determining the solution of
the absolute value equation is equal to determining the solution of
the linear system. That happens all the time, see the next section
(Examples 3.1 and 3.2). In such case, we can employ some numerical
methods, like Krylov subspace methods, for solving the corresponding
linear systems. That is to say, once we judge that the solution $x$
of the GAVE (\ref{eq:11}) is non positive or positive, i.e., $x\leq
0$ or $x> 0$,  we only needs to deal with the linear systems, no
longer deal with the GAVE (\ref{eq:11}).

Finally, we turn to discuss the convergence property of Method 3.
Theorem \ref{th1} presents a  general convergence condition of
Method 3 for solving the GAVE (\ref{eq:11}) when the related matrix
is nonsingular.

\begin{theorem} \label{th1}
Let $A+B+\Omega$ be invertible. Define
\[
f_{\Omega}=|(A+B+\Omega)^{-1}\Omega|,
g_{\Omega}=2|(A+B+\Omega)^{-1}B|.
\]
If
\begin{equation}\label{eq:23}
\rho(f_{\Omega}+g_{\Omega})<1,
\end{equation}
then Method 3 is convergent.
\end{theorem}
\textbf{Proof.} Assume that a solution of the GVAE (\ref{eq:11}) is
$x^{\ast}$. Combining (\ref{eq:21}) with (\ref{eq:22}), we get
\begin{equation}\label{eq34}
x^{k+1}-x^{\ast}=(A+B+\Omega)^{-1}(\Omega (x^{k}-x^{\ast})
+2B(\max\{0,x^{k}\}-\max\{0,x^{\ast}\}).
\end{equation}
Together with Corollary 2.1 in \cite{Li24}, using the absolute value
for (\ref{eq34}) leads to
\begin{align*}\label{eq:34}
|x^{k+1}-x^{\ast}|=&|(A+B+\Omega)^{-1}(\Omega (x^{k}-x^{\ast})
+2B(\max\{0,x^{k}\}-\max\{0,x^{\ast}\})|\\
\leq&|(A+B+\Omega)^{-1}\Omega (x^{k}-x^{\ast})|+2|(A+B+\Omega)^{-1}B(\max\{0,x^{k}\}-\max\{0,x^{\ast}\})|\\
\leq&|(A+B+\Omega)^{-1}\Omega||x^{k}-x^{\ast}|+2|(A+B+\Omega)^{-1}B||x^{k}-x^{\ast}|\\
=&(f_{\Omega}+g_{\Omega})|x^{k}-x^{\ast}|.
\end{align*}
Clearly, Method 3 is convergent under the condition (\ref{eq:23}).
$\hfill{} \Box$

Based on Theorem \ref{th1}, it is easy to see that Method 3 is also
convergent if $A,B$ and $\Omega$ satisfy
\[
\|(A+B+\Omega)^{-1}\|_{2}(\|\Omega\|_{2}+2\|B\|_{2})<1.
\]
Further, together with  Banach perturbation Lemma in \cite{Ortega},
Theorem \ref{th2} can be given.

\begin{theorem} \label{th2}
Let $A+B$ and $A+B+\Omega$ be invertible. If
\begin{equation}\label{eq:25}
\|(A+B)^{-1}\|_{2}<(2\|\Omega\|_{2}+2\|B\|_{2})^{-1},
\end{equation}
then Method 3 is convergent.
\end{theorem}
\textbf{Proof.} Based on  the Banach perturbation lemma in
\cite{Ortega}, we have
\begin{align*}
\|(A+B+\Omega)^{-1}\|_{2}&\leq\frac{\|(A+B)^{-1}\|_{2}}{1-
\|(A+B)^{-1}\|_{2}\|\Omega\|_{2}}<\frac{\frac{1}{2\|B\|_{2}+
2\|\Omega\|_{2}}}{1-\frac{\|\Omega\|_{2}}{2\|B\|_{2}+2\|\Omega\|_{2}}}=\frac{1}{\|\Omega\|_{2}+2\|B\|_{2}}.
\end{align*}
This shows that Method 3 is convergent under the condition
(\ref{eq:25}). $\hfill{} \Box$

In addition, we consider two cases: (1) matrix $A+B$ is symmetric
positive definite; (2) matrix $A+B$ is an $H_{+}$-matrix.

When matrix $A+B$ is symmetric positive definite, for $\Omega=\omega
I$ with $\omega>0$, Theorem \ref{th3} can be obtained.

\begin{theorem}\label{th3}
Let $A+B$ be  symmetric positive definite, $2\|B\|_{2}=\tau$ and
$\Omega=\omega I$ with $\omega>0$. Further, let $\mu_{\min}$ be the
smallest eigenvalue of the matrix $A+B$. If
\begin{equation}\label{eq:26}
\tau<\mu_{\min},
\end{equation}
then Method 3 is convergent.
\end{theorem}
\textbf{Proof.} By assumptions, we have
\begin{align*}
\|(A+B+\Omega)^{-1}\|_{2}(\|\Omega\|_{2}+2\|B\|_{2})&=\|(A+B+\omega
I)^{-1}\|_{2}(\|\omega
I\|_{2}+2\|B\|_{2})=\frac{\omega+\tau}{\mu_{\min}+\omega}.
\end{align*}
It is not difficult to check that  Method 3 is convergent under the
condition (\ref{eq:26}). $\hfill{} \Box$

The advantage of Theorem \ref{th3} is that it does not require that
$A$ and $B$ must be symmetric positive definite.

When matrix $A+B$ is an $H_{+}$-matrix, Theorem \ref{th4} can be
derived.

\begin{theorem}\label{th4}
Let $A+B$ be an $H_{+}$-matrix, and  $\Omega$ be a positive diagonal
matrix. If $\langle A+B\rangle-2|B|$ is an $M$-matrix, then Method 3
is convergent.
\end{theorem}
\textbf{Proof.} By assumptions and Lemma 2.1, we know
\[
|(A+B+\Omega)^{-1}|\leq(\langle A+B\rangle+\Omega)^{-1}.
\]
Based on the proof of Theorem \ref{th1}, we have
\begin{align*}
|x^{k+1}-x^{\ast}|&\leq|(A+B+\Omega)^{-1}|\cdot(\Omega
+2|B|)|x^{k}-x^{\ast} \leq(\langle A+B\rangle+\Omega)^{-1}(\Omega
+2|B|)|x^{k}-x^{\ast}|.
\end{align*}
Therefore, Method 3 is convergent if $\langle A+B\rangle-2|B|$ is an
$M$-matrix. $\hfill{} \Box$

Similar to Theorem \ref{th3}, the advantage of Theorem \ref{th4} is
that it does not require that $A$ and $B$ must be an $H_{+}$-matrix.
In particular, for $x^{\ast}$ being non positive (it is equal to the
nonlinear term vanished in Method 3), when  $A+B$ is an
$H_{+}$-matrix, Method 3 is convergent for any positive diagonal
matrix $\Omega$.

\section{Numerical experiments}
In this section, by using two examples, we show the performance of
Method 3 for solving the GAVE (1.1) from three aspects: the
iteration steps (indicated as `IT'), the elapsed CPU time (indicated
as `CPU') and the relative residual error (indicated as `RES').
Here, `RES' is set to be
\[
\mbox{RES}=\|Ax^{k}-B|x^{k}|-b\|_{2}/\|b\|_{2}.
\]

In our computations, the initial vector of the tested methods is set
to be zero vector and all iterations are stopped once
$\mbox{RES}<10^{-6}$ or the prescribed iteration count
$k_{\max}=500$ is surpassed. All the tests are operated on MATLAB
R2016b.

\textbf{Example 3.1} (\cite{Bai})  Consider the following GAVE
\begin{equation}\label{eq:31}
Ax-B|x|=q,
\end{equation}
where $A=R+I$, $B=R-I$ and $x=\frac{1}{2}((R-I)z+q)$, which is from
the linear complementarity problem (LCP)
\begin{equation}\label{eq:32}
z > 0, w = Rz + q > 0, z^{T}w = 0.
\end{equation}
Clearly, using $w=|x|+x$  and $z=|x|-x$ for the LCP (\ref{eq:32})
leads to the GAVE (\ref{eq:31}). In our computations, $R =
\hat{R}+\mu I$ with
\[
\hat{R}=\mbox{tridiag}(-I,S, -I)\in \mathbb{R}^{n\times n} \
\mbox{and}\ S=\mbox{tridiag}(-1,4, -1)\in \mathbb{R}^{m\times m},
\]
and $q = -Rz^{\ast}$ with $z^{\ast}=(1, 2, 1, 2,\ldots, 1,2)^{T}\in
\mathbb{R}^{n}$ being the unique solvability of the LCP($q,R$). In
this case, the unique solvability of the corresponding GAVE
(\ref{eq:31}) is $x^{\ast} = (-0.5,-1,\ldots, -0.5,-1)^{T} \in
\mathbb{R}^{n}$. Let $m$ be a prescribed positive integer and
$n=m^{2}$.

\textbf{Example 3.2} (\cite{Bai}) It is nonsymmetric version of
Example 3.1. That is to say, $R = \hat{R}+\mu I$  with
\[
\hat{R}=\mbox{tridiag}(-1.5I,S, -0.5I)\in \mathbb{R}^{n\times n} \
\mbox{and}\ S=\mbox{tridiag}(-1.5,4, -0.5)\in \mathbb{R}^{m\times
m},
\]
and $q = -Mz^{\ast}$ with $z^{\ast}=(1, 2, 1, 2,\ldots, 1, 2)^{T}\in
\mathbb{R}^{n}$.
%
%
and $q = -Mz^{\ast}$ with $z^{\ast}=(1, 2, 1, 2,\ldots, 1, 2)^{T}\in
\mathbb{R}^{n}$ being the unique solvability of the LCP($q,R$). In
this case, the unique solvability of the corresponding GAVE
(\ref{eq:31}) is $x^{\ast} = (-0.5,-1,\ldots, -0.5,-1)^{T} \in
\mathbb{R}^{n}$. Let $m$ be a prescribed positive integer and
$n=m^{2}$.

\begin{table}[!htb] \centering
\begin{tabular}
{p{60pt}p{40pt}p{50pt}p{50pt}p{50pt}p{50pt}} \hline
&$n$ &$2500$ &$10000$&$22500$&$40000$\\
 \hline M-1&  IT&17&18&18&18\\
&CPU&0.0106&  0.0735& 0.2305 &  0.5245\\
&RES&9.5311e-7&4.3408e-7& 4.4330e-7&4.4793e-7\\
 \hline
M-2&  IT&2&2&2&2  \\
&CPU&0.0029&0.0037&0.0042&0.0076\\
&RES&4.7468e-7& 2.1665e-7&2.2139e-7&2.2377e-7\\
 \hline
M-3& IT&2&2&2&2  \\
&CPU&0.0021&0.054&0.0140&0.0297\\
&RES&2.5077e-7& 1.1458e-7&1.1713e-7&1.1841e-7\\
 \hline
G-20& IT&10&10&10&10  \\
&CPU&0.0865&0.0836&0.0900&0.0952\\
&RES&8.4898e-7& 6.2799e-7&5.1978e-7&4.5312e-7\\
 \hline

\end{tabular}
\\ \caption{ Numerical comparison of Example 3.1 for $\mu=4$ with $\Omega=\mbox{diag}(A)$.}
\end{table}

For Examples 3.1 and 3.2, $x^{\ast}$ is non positive.
In such case, we take $y=x^{\ast}$, Eq. (\ref{eq:31}) reduces to
\begin{equation*}
(A+B+\Omega)y=\Omega y+q.
\end{equation*}
Further, we have the following method
\begin{equation*}
y^{k+1}=(A+B+\Omega)^{-1}(\Omega y^{k}+q),
\end{equation*}
which is equal to the nonlinear term vanished in Method 3. It is
easy to check that it is convergent because $A+B$ is  an
$H_{+}$-matrix. Moreover, we directly use Krylov subspace methods
like GMRES($m$) for the corresponding linear systems $(A+B)x=q$,
where `$m$' denotes the restated number of GMRES.

For Examples 3.1 and 3.2,  based on the different problem sizes of
$n$, we test Method 1, the Gauss-Seidel version of Method 2, Method
3 and GMRES(20) for $\mu=4$. In the tables, `M-1', `M-2', `M-3' and
`G-20' in order denote Method 1, the Gauss-Seidel version of Method
2, Method 3 and GMRES(20).

\begin{table}[!htb] \centering
\begin{tabular}
{p{60pt}p{40pt}p{50pt}p{50pt}p{50pt}p{50pt}} \hline
&$n$ &$2500$ &$10000$&$22500$&$40000$\\
 \hline M-1&  IT&17 &18 &18 &18 \\
&CPU&0.0105 &0.0733 &0.2328 &0.5178 \\
&RES&9.3131e-7 & 4.2889e-7 &4.3983e-7 &4.4533e-7 \\
 \hline
M-2&  IT&2 &2 &2 &2 \\
&CPU&0.0021 &0.0024 & 0.0040 &0.0060 \\
&RES&4.3163e-7 &1.9951e-7 &2.0483e-7 &2.0750e-7 \\
 \hline
M-3& IT&2 &2 &2 &2 \\
&CPU&0.0024 &0.0053 &0.0139 &0.0299
 \\
&RES&2.2775e-7 &1.0546e-7 &1.0833e-7 &1.0977e-7 \\
 \hline
G-20 & IT&12&12&12&12  \\
&CPU&0.0846&0.0884&0.1008&0.1080\\
&RES&8.9460e-7&  6.1218e-7&4.9215e-7&4.2247e-7\\
 \hline
\end{tabular}
\\ \caption{ Numerical comparison of Example 3.2 for $\mu=4$ with $\Omega=\mbox{diag}(A)$.}
\end{table}

Here, for the former three  methods, the choice of $\Omega$ is
$\Omega=\mbox{diag}(A)$. Using the former three methods, we need to
deal with the linear systems with the matrices $A$, $D+\Omega-L$ and
$A+B+\Omega$, respectively. In the implementations, for solving
these linear systems, the sparse LU factorization can be employed.

Tables 1 and 2 list the numerical results  (including IT, CPU and
RES) of all the tested methods. From these numerical results,  we
observe that all tested methods can quickly compute a satisfactory
approximation to the solution of the corresponding GAVE
(\ref{eq:31}). Further, observing these numerical results in  Tables
1 and 2, we have the following facts:

\begin{itemize}
\item Method 3 and  the Gauss-Seidel version of Method 2 require less
iteration steps and computing times than Method 1. This shows that
Method 3 and  the Gauss-Seidel version of Method 2 outperforms
Method 1.
\item Comparing Method 3 and the Gauss-Seidel version of Method 2, both requires the same
iteration steps. Although the computing times of the former are a
little more than the latter, the relative residual error of the
former is less than the latter. This shows the former is more
accurate than the latter. Roughly speaking, Method 3 is equal to the
Gauss-Seidel version of  Method 2 in a way (in fact, this version of
Method 2 is its best relaxation version, see \cite{Zhou21}).
\item Comparing Method 3 with GMRES(20), it is easy to find that
Method 3 overmatches  GMRES(20) in terms of the iteration steps, the
computing times and the relative residual error. In terms of
computing efficiency Method 3 is superior to GMRES(20).
\end{itemize}

In total, we can see that Method 3 for the GAVE (\ref{eq:11}) is
with good performance, that is to say, it is feasible and
competitive, compared with some existing iteration methods.

\section{Conclusions}
In this paper, a class of maximum-based iteration methods has been
designed to solve the GAVE (\ref{eq:11}) by using
$|x|=2\max\{0,x\}-x$. Its convergence conditions are presented. What
we find interesting is that under certain condition determining the
solution of the absolute value equation is equal to determining the
solution of the linear system. Numerical experiments verify that
this new iteration method are superior to some existing iteration
methods, including the modified Newton method, some Krylov subspace
methods like GMRES(20).



\end{document}